\documentclass[11pt]{article}
\usepackage{amsmath}
\usepackage{amssymb}
\usepackage{amsthm}
\usepackage{hyperref}
\usepackage[german,english]{babel}

\newcommand\nonu{\nonumber}
\newcommand\bPP{\\[\bigskipamount]\indent}
\newcommand\RR{\mathbb{R}}
\newcommand\al\alpha
\newcommand\be\beta
\newcommand\ga\gamma
\newcommand\de\delta
\newcommand\la\lambda
\newcommand\Ga{\Gamma}
\newcommand\tha\theta
\newcommand\half{\frac12}
\newcommand\thalf{\tfrac12}
\newcommand\iy\infty
\newcommand\const{\mathrm{const.}\,}
\newcommand\LHS{left-hand side}
\newcommand\RHS{right-hand side}
\renewcommand\Re{\operatorname{Re}}
\renewcommand\Im{\operatorname{Im}}
\newcommand{\hyp}[5]{\,\mbox{}_{#1}F_{#2}\!\left(
  \genfrac{}{}{0pt}{}{#3}{#4};#5\right)}

\makeatletter
\newcommand{\dlmf}[1]{%
\cite[%
 \def\nextitem{\def\nextitem{, }}%
 \@for \el:=#1\do{\nextitem\expandafter\dlmf@eq@href\el...\end}%
]{dlmf}%
}
\def\dlmf@eq@href#1.#2.#3.#4\end{%
  \href{http://dlmf.nist.gov/#1.#2.E#3}{(#1.#2.#3)}}
\makeatother

\numberwithin{equation}{section}
\newtheorem{theorem}{Theorem}[section]

\newtheorem{Definition}[theorem]{Definition}

\newtheorem{Remark}[theorem]{Remark}
\newenvironment{remark}{\begin{Remark}\rm}{\end{Remark}}
\newtheorem{Example}[theorem]{Example}

\newcommand\Proof{\noindent{\bf Proof}\quad}

\begin{document}
\title{Dual addition formulas associated with dual product formulas}
\author{Tom H. Koornwinder}
\date{Dedicated to Mourad E. H. Ismail on the occasion of his
seventieth birthday}
\maketitle
\begin{abstract}
We observe that the linearization coefficients for ultraspherical polynomials
are the orthogonality weights for Racah polynomials with special parameters. Then it turns out that
the linearization sum with such a Racah polynomial as extra factor inserted,
can also be evaluated. The corresponding Fourier--Racah expansion is
an addition type formula which is dual to the well-known addition
formula for ultraspherical polynomials. The limit to the case of
Hermite polynomials of this dual addition formula
is also considered.
Similar results as for ultraspherical polynomials, although only formal, 
are given by taking the Ruijsenaars--Halln\"as dual product formula for
Gegenbauer functions as a starting point and by working with Wilson
polynomials. 
\end{abstract}
\section{Introduction}
A prototype for an addition formula for a family of special
orthogonal polynomials is the
\emph{ addition formula for Legendre polynomials}
\dlmf{18.18.9}
\begin{multline}
P_n(\cos\tha_1\cos\tha_2+\sin\tha_1\sin\tha_2\cos\phi)
=P_n(\cos\tha_1)P_n(\cos\tha_2)\\
+2\sum_{k=1}^n \frac{(n-k)!\,(n+k)!}{2^{2k}(n!)^2}\,
(\sin\tha_1)^k P_{n-k}^{(k,k)}(\cos\tha_1)
\,(\sin\tha_2)^k P_{n-k}^{(k,k)}(\cos\tha_2)\cos(k\phi).
\label{35}
\end{multline}
The \RHS\ is the Fourier-cosine expansion of the \LHS\ as a function
of $\phi$. Integration with respect to $\phi$ over $[0,\pi]$
gives the constant term in this expansion, \dlmf{18.17.6}
\begin{equation}
P_n(\cos\tha_1)P_n(\cos\tha_2)=\frac1\pi \int_0^\pi
P_n(\cos\tha_1\cos\tha_2+\sin\tha_1\sin\tha_2\cos\phi)\,d\phi,
\label{36}
\end{equation}
which is known as the \emph{ product formula for Legendre polynomials}.

A formula dual to \eqref{36} is the \emph{ linearization formula
for Legendre polynomials}, see \dlmf{18.18.22} for $\la=\half$
together with \dlmf{18.7.9} and see \dlmf{5.2.4} for the shifted factorial $(a)_n$. It reads:
\begin{equation}
P_l(x)P_m(x)=\sum_{j=0}^{\min(l,m)}
\frac{(\half)_j(\half)_{l-j}(\half)_{m-j} (l+m-j)!}
{j!\,(l-j)!\,(m-j)!\,(\tfrac32)_{l+m-j}}\,\big(2(l+m-2j)+1\big)\,
P_{l+m-2j}(x).
\label{37}
\end{equation}
On several occasions, during his lectures at conferences, Richard Askey
raised the problem to find an addition type formula associated with
\eqref{37} in a similar way as the addition formula \eqref{35} is
associated with the product formula \eqref{36}. It is the purpose of the
present paper to give such a formula, more generally associated with
the linearization formula for ultraspherical polynomials, and also a
(formal) addition type formula associated with the dual product
formula for Gegenbauer functions which was recently given by
Halln\"as \& Ruijsenaars \cite[(4.17)]{1}.

In order to get a better feeling for what a dual addition formula should
look like, we first rewrite \eqref{35} and \eqref{36} by substituting
$z=\sin\tha_1\sin\tha_2\cos\phi$, $x=\cos\tha_1$,
$y=\cos\tha_2$, and putting $\cos(k\phi)=T_k(\cos\phi)$
($T_k$ is a Chebyshev polynomial \dlmf{18.5.1}). Assume that
$x,y\in[-1,1]$ and
$z\in\big[-\sqrt{1-x^2}\,\sqrt{1-y^2},\sqrt{1-x^2}\,\sqrt{1-y^2}\,\big]$.
We obtain
\begin{multline}
P_n(z+xy)=P_n(x)P_n(y)
+2\sum_{k=1}^n \frac{(n-k)!\,(n+k)!}{2^{2k}(n!)^2}\\
\times (1-x^2)^{\half k} P_{n-k}^{(k,k)}(x)\,
(1-y^2)^{\half k} P_{n-k}^{(k,k)}(y)\,
T_k\left(\frac z{\sqrt{1-x^2}\,\sqrt{1-y^2}}\right)
\label{38}
\end{multline}
and
\begin{equation}
P_n(x)P_n(y)=\frac1\pi
\int_{-\sqrt{1-x^2}\,\sqrt{1-y^2}}^{\sqrt{1-x^2}\,\sqrt{1-y^2}}\,
\frac{P_n(z+xy)}{\sqrt{(1-x^2)(1-y^2)-z^2}}\,dz.
\label{39}
\end{equation}
The rewritten addition formula \eqref{38} expands the \LHS\
as a function of $z$
in terms of Chebyshev polynomials of dilated argument, where the dilation
factor depends on $x,y$, and the product formula \eqref{39}
recovers the constant term in this expansion by integration with
respect to the weight function over the orthogonality interval for
these dilated Chebyshev polynomials.
The linearization formula \eqref{37} very much looks as a dual formula
with respect to \eqref{39}. If we can recognize the coefficents in the sum
in \eqref{37} as weights for some finite system of orthogonal polynomials
then the dual addition formula associated with \eqref{37} should be the
corresponding orthogonal expansion of $P_{l+m-2j}(x)$ as a function of $j$.

It is not so easy to recognize the coefficients in \eqref{37}
as weights for known orthogonal polynomials, but a strong hint was
provided by the Halln\"as-Ruijsenaars dual product formula
\cite[(4.17)]{1} for Gegenbauer functions, which can be rewritten as
\eqref{7}. There the weight function in the integral is clearly
the weight function for Wilson polynomials \cite[Section 9.1]{9}
with suitable parameters.
This suggests that in \eqref{37} we should have weights of
Racah polynomials \cite[Section~9.2]{9}, whch are the discrete analogues
of the Wilson polynomials. Indeed, this turns out to work, see \eqref{17}
(more generally for ultraspherical polynomials),
and a nice expansion \eqref{40}
in terms of these Racah polynomials can be derived,
which is the dual addition formula for ultraspherical polynomials
predicted by Askey. In a remark at the end of Section \ref{61}
we point to a paper by Koelink et al.\ \cite{17} from 2013
which already has the dual addition formula in disguised form.

For better comparison of the dual results in Section \ref{61}
we state the addition formula for ultraspherical polynomials and
related formulas in Section \ref{60}. There we also mention
the quite unknown paper \cite{18} by All\'e from 1865 which already
gives the addition formula for ultraspherical polynomials, much earlier
than Gegenbauer's paper \cite{19} from 1874.

In Section \ref{62} we obtain a limit of the dual addition formula
for ultraspherical polynomials corresponding to the limit from
ultraspherical to Hermite polynomials.
The resulting formulas for Hermite polynomials are well known.
Remarkable is that the orthogonality for special Racah polynomials
tends in the limit to a (well known)
biorthogonality for shifted factorials.

In Section \ref{63} we introduce the needed special functions.
The paper concludes in Section \ref{64} with a list of possible
follow-up work on dual addition formulas.
\section{Preliminaries}
\label{63}
\subsection{Jacobi, ultraspherical and Hermite polynomials}
We will work with renormalized \emph{ Jacobi polynomials}
\cite[Section 9.8]{9}
\begin{equation}
R_n^{(\al,\be)}(x):=\frac{P_n^{(\al,\be)}(x)}{P_n^{(\al,\be)}(1)}
=\hyp21{-n,n+\al+\be+1}{\al+1}{\thalf(1-x)},
\label{27}
\end{equation}
where $P_n^{(\al,\be)}(1)=(\al+1)_n/n!$ and $R_n^{(\al,\be)}(1)=1$
(see \dlmf{16.2.1} for the definition of a ${}_pF_q$ hypergeometric
series).
These are orthogonal
polynomials on the interval $[-1,1]$ with respect to the weight function
$(1-x)^\al(1+x)^\be$ ($\al,\be>-1$). In particular, for $\al=\be$,
the polynomials are called \emph{ ultraspherical} or
\emph{Gegenbauer polynomials}, then the weight
function is even and we have
\begin{equation*}
R_n^{(\al,\al)}(-x)=(-1)^n R_n^{(\al,\al)}(x).
\end{equation*}
The precise orthogonaliy relation is
\begin{equation}
\begin{split}
&\int_{-1}^1 R_m^{(\al,\al)}(x)\,R_n^{(\al,\al)}(x)\,(1-x^2)^\al\,dx
=h_n^{(\al,\al)}\,\de_{m,n}\,,\\
&\qquad\qquad\qquad\qquad
h_n^{(\al,\al)}=\frac{2^{2\al+1}\Ga(\al+1)^2}{\Ga(2\al+2)}\,
\frac{n+2\al+1}{2n+2\al+1}\,\frac{n!}{(2\al+2)_n}\,.
\end{split}
\label{57}
\end{equation}
From \eqref{27} we see that
\begin{equation}
R_n^{(\al,\be)}(x)=\frac{(n+\al+\be+1)_n}{2^n(\al+1)_n}\,x^n+
\mbox{terms of lower degree}.
\label{28}
\end{equation}
The connection of $R_n^{(\al,\al)}$ with the usual $C_n^{(\la)}$
notation for ultraspherical polynomials is:
\begin{equation*}
R_n^{(\al,\al)}(x)=\frac{n!}{(\al+1)_n}\,P_n^{(\al,\al)}(x)
=\frac{n!}{(2\al+1)_n}\,C_n^{(\al+\half)}(x).
\end{equation*}
From \dlmf{18.5.10} we have the power series
\begin{equation}
R_n^{(\al,\al)}(x)=\frac{n!}{(2\al+1)_n}
\sum_{k=0}^{[\half n]} \frac{(-1)^k (\al+\thalf)_{n-k}}
{k!\,(n-2k)!}\,(2x)^{n-2k}
\label{50}
\end{equation}

We will need the difference formula
\begin{equation}
R_n^{(\al,\al)}(x)-R_{n-2}^{(\al,\al)}(x)
=\frac{n+\al-\half}{\al+1}\,(x^2-1)\,R_{n-2}^{(\al+1,\al+1)}(x)\qquad
(n\ge2).
\label{23}
\end{equation}
{\bf Proof of \eqref{23}.}\quad
More generally, let $w(x)=w(-x)$ be an even weight function on $[-1,1]$,
let $p_n(x)=k_nx^n+\cdots\;$ be orthogonal polynomials on $[-1,1]$ with
respect to the weight function $w(x)$, and let
$q_n(x)=k_n'x^n+\cdots\;$ be orthogonal polynomials on $[-1,1]$ with
respect to the weight function $w(x)(1-x^2)$. Assume that $p_n$ and $q_n$
are normalized by $p_n(1)=1=q_n(1)$. Let $n\ge2$.
Then $p_n(x)-p_{n-2}(x)$ vanishes for $x=\pm1$. Hence
$(p_n(x)-p_{n-2}(x))/(1-x^2)$ is a polynomial of degree $n-2$. It is seen
immediately that $x^k$ ($k<n-2$) is orthogonal to this polynomial with
respect to the weight function $w(x)(1-x^2)$ on $[-1,1]$. We conclude that
\[
p_n(x)-p_{n-2}(x)=\frac{k_n}{k_{n-2}'}\,(x^2-1)q_{n-2}(x)\qquad(n\ge2).
\]
Now specialize to $w(x)=(1-x^2)^\al$ and use \eqref{28}.\qed

\emph{ Hermite polynomials} \cite[Section 9.15]{9} are orthogonal polynomials
$H_n$ on $(-\iy,\iy)$ with respect to the weight function $e^{-x^2}$
and normalized such that $H_n(x)=2^nx^n+\mbox{terms of lower degree}$.
From \dlmf{18.5.13} we have the power series
\begin{equation}
H_n(x)=n!\,\sum_{k=0}^n \frac{(-1)^k}{k!\,(n-2k)!}\,(2x)^{n-2k}
\label{51}
\end{equation}
It follows from \eqref{50} and \eqref{51} that
\begin{align}
\lim_{\al\to\iy}\al^{\half n} R_n^{(\al,\al)}(\al^{-\half}x)
&=2^{-n} H_n(x),
\label{52}\\
\lim_{\al\to\iy}\al^{\mu n} R_n^{(\al,\al)}(\al^{-\mu} x)&=x^n
\qquad(\mu<\thalf),
\nonu\\
\noalign{\noindent\mbox{in particular,}}
\lim_{\al\to\iy}R_n^{(\al,\al)}(x)&=x^n.
\label{53}
\end{align}
\subsection{Racah polynomials}
We will consider \emph{ Racah polynomials}
\cite[Section 9.2]{9}
\begin{equation}
R_n\big(x(x+\ga+\de+1);\al,\be,\ga,\de\big):=
\hyp43{-n,n+\al+\be+1,-x,x+\ga+\de+1}{\al+1,\be+\de+1,\ga+1}1
\label{24}
\end{equation}
for $\ga=-N-1$, where $N\in\{1,2,\ldots\}$, and for
$n\in\{0,1,\ldots,N\}$. These are orthogonal polynomials on
the finite quadratic set $\{x(x+\ga+\de+1)\mid x\in\{0,1,\ldots,N\}\}$:
\begin{equation}
\sum_{x=0}^N (R_mR_n)\big(x(x+\ga+\de+1);\al,\be,\ga,\de\big)\,
w_{\al,\be,\ga,\de}(x)=h_{n;\al,\be,\ga,\de} \de_{m,n}\quad
(m,n\in\{0,1,\ldots,N\})
\label{30}
\end{equation}
with
\begin{equation}
w_{\al,\be,\ga,\de}(x)=
\frac{(\al+1)_x (\be+\de+1)_x (\ga+1)_x (\ga+\de+1)_x}
{(-\al+\ga+\de+1)_x (-\be+\ga+1)_x (\de+1)_x\,x!}\,
\frac{\ga+\de+1+2x}{\ga+\de+1},
\label{26}
\end{equation}
\begin{equation}
\frac{h_{n;\,\al,\be,\ga,\de}}{h_{0;\,\al,\be,\ga,\de}}=
\frac{\al+\be+1}{\al+\be+2n+1}\,
\frac{(\be+1)_n(\al+\be-\ga+1)_n(\al-\de+1)_n\,n!}
{(\al+1)_n(\al+\be+1)_n(\be+\de+1)_n(\ga+1)_n},
\label{31}
\end{equation}
\begin{equation}
h_{0;\,\al,\be,\ga,\de}=\sum_{x=0}^N w_{\al,\be,\ga,\de}(x)=
\frac{(\al+\be+2)_N(-\de)_N}{(\al-\de+1)_N(\be+1)_N}\,.
\label{29}
\end{equation}

Clearly $R_n(0;\al,\be,\ga,\de)=1$ while, by \eqref{24} and
the Saalsch\"utz formula \dlmf{16.4.3}, we can evaluate
the Racah polynomial for $x=N$:
\begin{equation}
R_n(N\de;\al,\be,\ga,\de)=
\frac{(\be+1)_n(\al-\de+1)_n}{(\al+1)_n(\be+\de+1)_n}\,.
\label{25}
\end{equation}

The backward shift operator equation
\cite[(9.2.8)]{9} can be rewritten as
\begin{align}
&w_{\al,\be,\ga,\de}(x) R_n\big(x(x+\ga+\de+1);\al,\be,\ga,\de\big)\nonu\\
&=\frac{\ga+\de+2}{\ga+\de+2+2x}\,w_{\al+1,\be+1,\ga+1,\de}(x)\,
R_{n-1}\big(x(x+\ga+\de+2);\al+1,\be+1,\ga+1,\de\big)\nonu\\
&-\frac{\ga+\de+2}{\ga+\de+2x}\,w_{\al+1,\be+1,\ga+1,\de}(x-1)\,
R_{n-1}\big((x-1)(x+\ga+\de+1);\al+1,\be+1,\ga+1,\de\big).
\label{20}
\end{align}
This holds for $x=0,\ldots,N$. For $x=0$ \eqref{20} remains true if
we put the second term on the right equal to 0, while for $x=N$
the first term on the right can be assumed to vanish.
In this last case the identity \eqref{20} can be checked by using
\eqref{25} and \eqref{26}.

Hence, for a function $f$ on $\{0,1,\ldots,N\}$ we have
\begin{multline}
\sum_{x=0}^N
w_{\al,\be,\ga,\de}(x) R_n\big(x(x+\ga+\de+1);\al,\be,\ga,\de\big)\,f(x)
=\sum_{x=0}^{N-1}
\frac{\ga+\de+2}{\ga+\de+2+2x}\\
\times w_{\al+1,\be+1,\ga+1,\de}(x)\,
R_{n-1}\big(x(x+\ga+\de+2);\al+1,\be+1,\ga+1,\de\big)\,
\big(f(x)-f(x+1)\big).
\label{21}
\end{multline}
\subsection{Jacobi and Gegenbauer functions}
In \cite[(4.3)]{1} Halln\"as \& Ruijsenaars define a
\emph{ conical function}
\begin{equation}
F(g;r,2k):=
\left(\frac\pi 4\right)^\half
\frac{\Ga(g+ik)\Ga(g-ik)}{\Ga(g) (2\sinh r)^{g-\half}}\,
P_{ik-\half}^{\half-g}(\cosh r)\qquad
(r>0,\;\Re g>0).
\label{1}
\end{equation}
Here the $P$-function is the \emph{ associated Legendre function
of the first kind} which is expressed by \dlmf{14.3.6}
and \dlmf{15.1.2} as
Gauss hypergeometric function:
\begin{equation}
P_\nu^\mu(x)=\frac1{\Ga(1-\mu)}\left(\frac{x+1}{x-1}\right)^{\half\mu}
\hyp21{\nu+1,-\nu}{1-\mu}{\thalf-\thalf x}\qquad(x>1).
\label{2}
\end{equation}
Substitute \eqref{2} in \eqref{1} and also use Euler's
transformation formula \dlmf{15.8.1}
and Legendre's duplication formula \dlmf{5.5.5}.
Then we obtain
\begin{equation}
F(g;r,2k)=\frac{\Ga(g+ik)\Ga(g-ik)}{2\Ga(2g)}\,
\hyp21{g+ik,g-ik}{g+\thalf}{-\sinh^2\thalf r}.
\end{equation}
By \dlmf{15.9.11} this can be written in terms of \emph{ Jacobi functions}
\cite{2}, \cite{3}
\begin{equation}
\phi_\la^{(\al,\be)}(t):=
\hyp21{\thalf(\al+\be+1+i\la),\thalf(\al+\be+1-i\la)}{\al+1}{-\sinh^2 t}
\qquad(t\in\RR)
\label{11}
\end{equation}
(called \emph{ Gegenbauer functions} if $\al=\be$) as
\begin{equation*}
F(g;r,2k)=\frac{\Ga(g+ik)\Ga(g-ik)}{2\Ga(2g)}\,
\phi_{2k}^{(g-\half,g-\half)}(\thalf r).
\end{equation*}
By \cite[(2.8)]{2}
\begin{equation}
\phi_{2\la}^{(\al,\al)}(t)=\phi_\la^{(\al,-\half)}(2t),
\label{16}
\end{equation}
this becomes
\begin{equation}
F(g;r,2k)=\frac{\Ga(g+ik)\Ga(g-ik)}{2\Ga(2g)}\,
\phi_k^{(g-\half,-\half)}(r)
\label{4}
\end{equation}
or equivalently,
\begin{equation}
F(\al+\thalf;t,2\la)=
\frac{\Ga(\al+\thalf+i\la)\Ga(\al+\thalf-i\la)}{2\Ga(2\al+1)}\,
\phi_\la^{(\al,-\half)}(t).
\end{equation}
Note that, by \eqref{11}, $\phi_\la^{(\al,\be)}(0)=1$.
From \cite[(6.1)]{3} we have
\begin{equation}
|\phi_\la^{(\al,\be)}(t)|\le 1\qquad
(\al\ge\be\ge-\thalf,\;t\in\RR,\;|\Im\la|\le\al+\be+1).
\label{32}
\end{equation}

The contiguous relation
\begin{equation*}
\hyp21{a,b}cz-\hyp21{a-1,b+1}cz=\frac{(b-a+1)z}c\,\hyp21{a,b+1}{c+1}z,
\end{equation*}
which follows immediately by substitution of the power series
for the three ${}_2F_1$ functions, can be rewritten by \eqref{11}
in terms of Jacobi functions:
\begin{equation}
\frac{\phi_{\la-i}^{(\al,\be)}(t)-\phi_{\la+i}^{(\al,\be)}(t)}{i\la}=
\frac{\sinh^2 t}{\al+1}\,\phi_\la^{(\al+1,\be)}(t).
\label{34}
\end{equation}
\subsection{Wilson polynomials}
\emph{ Wilson polynomials} \cite[Section 9.1]{9} are defined by
\begin{equation}
\frac{W_n(x^2;a,b,c,d)}{(a+b)_n(a+c)_n(a+d)_n}:=
\hyp43{-n,n+a+b+c+d-1,a+ix,a-ix}{a+b,a+c,a+d}1.
\label{14}
\end{equation}
We need these polynomials
with parameters $\pm i\la\pm i\mu+\thalf\al+\tfrac14$
($\al>-\half$, $\la,\mu\in\RR$).
Then the orthogonality relation becomes \cite[(9.1.2)]{9}
\begin{multline}
\frac1{4\pi}\int_{-\iy}^\iy
(W_mW_n)(\nu^2;\pm i\la\pm i\mu+\thalf\al+\tfrac14)\,
\Bigg|\frac{\Ga\big(i\nu\pm i\la\pm i\mu+\thalf\al+\tfrac14\big)}
{\Ga(2i\nu)}\Bigg|^2 d\nu\\
=\frac{\Ga(\al+\thalf)^2\,|\Ga(n+\al+\thalf+2i\la)|^2\,
|\Ga(n+\al+\thalf+2i\mu)|^2}{\Ga(2n+2\al+1)}\,(n+2\al)_n\,n!\,\de_{m,n}.
\label{8}
\end{multline}
Here and later $\pm i\la\pm i\mu+\thalf\al+\tfrac14$ in the parameter
list means four elements in the list with the four possibilities
given by the two $\pm$ signs. Similarly
$\Ga\big(i\nu\pm i\la\pm i\mu+\thalf\al+\tfrac14\big)$
stands for a product of four Gamma functions.
We also wrote $(4\pi)^{-1}\int_{-\iy}^\iy$ instead of the usual
$(2\pi)^{-1}\int_0^\iy$, which is allowed because the integrand is
an even function of $\nu$.

The backward shift operator equation \cite[(9.1.9)]{9}
can be rewritten as
\begin{align}
&\frac{\Ga(a\pm ix)\ldots\Ga(d\pm ix)}{\Ga(\pm 2ix)}\,W_n(x^2;a,b,c,d)
\nonu\\
&\quad
=\frac{\Ga(a+\thalf\pm i(x+\thalf i))\ldots\Ga(d+\thalf+\pm i(x+\thalf i))}
{2i(x+\thalf i)\,\Ga(\pm 2i(x+\thalf i))}\,
W_n((x+\thalf i)^2;a+\thalf,b+\thalf,c+\thalf,d+\thalf)\nonu\\
&\quad
-\frac{\Ga(a+\thalf\pm i(x-\thalf i))\ldots\Ga(d+\thalf+\pm i(x-\thalf i))}
{2i(x-\thalf i)\,\Ga(\pm 2i(x-\thalf i))}\,
W_n((x-\thalf i)^2;a+\thalf,b+\thalf,c+\thalf,d+\thalf).
\label{33}
\end{align}
\section{The addition formula for ultraspherical polynomials}
\label{60}
In this section we briefly review the addition formula for ultraspherical
polynomials and formulas associated with it. This extends the discussion
of the Legendre case in the Introduction. As a reference see for instance
\cite[Section 9.8]{12}. We use the notation \eqref{27}.
\paragraph{Product formula} ($\al>-\half$)
\begin{equation}
R_n^{(\al,\al)}(x)R_n^{(\al,\al)}(y)=
\frac{\Ga(\al+1)}{\Ga(\al+\half)\Ga(\half)}
\int_{-1}^1 R_n^{(\al,\al)}\big(xy+(1-x^2)^\half(1-y^2)^\half t\big)\,
(1-t^2)^{\al-\half}\,dt.
\label{41}
\end{equation}
\paragraph{Addition formula}
\begin{multline}
R_n^{(\al,\al)}\big(xy+(1-x^2)^\half(1-y^2)^\half t\big)
=\sum_{k=0}^n\binom nk \frac{\al+k}{\al+\half k}\,
\frac{(n+2\al+1)_k (2\al+1)_k}{2^{2k} (\al+1)_k^2}\\
\times(1-x^2)^{\half k}R_{n-k}^{(\al+k,\al+k)}(x)\,
(1-y^2)^{\half k}R_{n-k}^{(\al+k,\al+k)}(y)\,R_k^{(\al-\half,\al-\half)}(t).
\label{42}
\end{multline}
\paragraph{Addition formula for $t=1$}
\begin{multline}
R_n^{(\al,\al)}\big(xy+(1-x^2)^\half(1-y^2)^\half\big)
=\sum_{k=0}^n\binom nk \frac{\al+k}{\al+\half k}\,
\frac{(n+2\al+1)_k (2\al+1)_k}{2^{2k} (\al+1)_k^2}\\
\times(1-x^2)^{\half k}R_{n-k}^{(\al+k,\al+k)}(x)\,
(1-y^2)^{\half k}R_{n-k}^{(\al+k,\al+k)}(y).
\label{44}
\end{multline}
For $x=\cos\tha_1$, $y=\cos\tha_2$ the \LHS\ takes the form
$R_n^{(\al,\al)}\big(\cos(\tha_1-\tha_2)\big)$.
\paragraph{Addition formula for $t=1$, $x=y$}
\begin{equation}
1=\sum_{k=0}^n\binom nk \frac{\al+k}{\al+\half k}\,
\frac{(n+2\al+1)_k (2\al+1)_k}{2^{2k} (\al+1)_k^2}\,
(1-x^2)^k \big(R_{n-k}^{(\al+k,\al+k)}(x)\big)^2.
\label{49}
\end{equation}
This shows in particular that $|R_n^{(\al,\al)}(x)|\le 1$ if $x\in[-1,1]$
and $\al>-\half$, \dlmf{18.14.1}. This is also well-known by
several other methods, including as a corollary of \eqref{41}.

\paragraph{Limit to Hermite polynomials}\quad\\
In the addition formula \eqref{42}
replace $x$ by $\al^{-\half}x$,
$t$ by $\al^{-\half} t$, multiply both sides of \eqref{42} by
$\al^{\half n}$ and let $\al\to\iy$. By \eqref{52} and \eqref{53}
we obtain
\begin{equation}
H_n\big(xy+(1-y^2)^\half t\big)=\sum_{k=0}^n \binom nk H_{n-k}(x)\,H_k(t)\,
(1-y^2)^{\half k}\,y^{n-k}.
\end{equation}
This is the case $n=2$ of \cite[10.13(40)]{16}.
The corresponding limit of the product formula \eqref{41} is
\begin{equation}
H_n(x)\,y^n=\frac1{\sqrt\pi}\int_{-\iy}^\iy H_n\big(xy+(1-y^2)^\half t\big)\,
e^{-t^2}\,dt.
\end{equation}

\paragraph{History of the addition formula for ultraspherical
polynomials}\quad\\
The addition formula \eqref{42} for ultraspherical polynomials is
usually ascribed to Gegenbauer \cite{19} (1874). However, it is
already stated and proved by All\'e \cite{18} in 1865.
The subsequent proofs by Gegenbauer \cite{19}, \cite{20}
in 1874 and 1893, and by Heine \cite[p.~455]{21} in (1878) do not
mention All\'e's result.
\section{The dual addition formula for ultraspherical polynomials}
\label{61}
The linearization formula for ultraspherical polynomials,
see \cite[(5.7)]{10}, can be written as
\begin{multline}
R_l^{(\al,\al)}(x) R_m^{(\al,\al)}(x)=
\frac{l!\,m!}{(2\al+1)_l(2\al+1)_m}
\sum_{j=0}^{\min(l,m)}\frac{l+m+\al+\half-2j}{\al+\half}\\
\times
\frac{(\al+\half)_j(\al+\half)_{l-j}(\al+\half)_{m-j}(2\al+1)_{l+m-j}}
{j!\,(l-j)!\,(m-j)!\,(\al+\frac32)_{l+m-j}}\,
R_{l+m-2j}^{(\al,\al)}(x).
\label{18}
\end{multline}
As mentioned in \cite[(4.18)]{11},
Rogers already gave the analogous linearization formula for
$q$-ultraspherical polynomials in 1895 and observed \eqref{18}
as a special case. Then \eqref{18} was independently given by Dougall
in 1919 without proof. See \cite[p.40]{10} for a discussion of further
treatments of \eqref{18}. See also \cite[Theorem 6.8.2]{12}
and the proof and discussion following the theorem.

From now on assume $\al>-\half$. Then the linearization coefficients
in \eqref{18} are nonnegative (as they are in the degenerate case
$\al=-\half$).
We also assume, without loss of generality, that $l\ge m$.

It is rather hidden in \eqref{18} that the linearization coefficients
are special cases of orthogonality weights \eqref{26} for
Racah polynomials. But indeed, 
a further rewriting of \eqref{18} and substitution of \eqref{26} and
\eqref{29} gives:
\begin{equation}
R_l^{(\al,\al)}(x) R_m^{(\al,\al)}(x)=
\sum_{j=0}^m
\frac{w_{\al-\half,\al-\half,-m-1,-l-\al-\half}(j)}
{h_{0;\,\al-\half,\al-\half,-m-1,-l-\al-\half}}\,
R_{l+m-2j}^{(\al,\al)}(x)\quad(l\ge m).
\label{17}
\end{equation}
This identity can be considered as giving the constant term
of an expansion of $R_{l+m-2j}^{(\al,\al)}(x)$ as a function of $j$
in terms of Racah polynomials
\begin{equation}
R_n\big(j(j-l-m-\al-\thalf);\,
\al-\thalf,\al-\thalf,-m-1,-l-\al-\thalf\big)
=\hyp43{-n,n+2\al,-j,j-l-m-\al-\thalf}{\al+\thalf,-l,-m}1.
\label{59}
\end{equation}
The general terms of this expansion will be obtained by evaluating
the sum
\begin{multline}
S_{n,l,m}^{(\al)}(x):=
\sum_{j=0}^m
w_{\al-\half,\al-\half,-m-1,-l-\al-\half}(j)\,
R_{l+m-2j}^{(\al,\al)}(x)\\
\times R_n\big(j(j-l-m-\al-\thalf);\,
\al-\thalf,\al-\thalf,-m-1,-l-\al-\thalf\big),
\label{22}
\end{multline}
where we still assume $l\ge m$ and where $n\in\{0,\ldots,m\}$.
\begin{theorem}
The sum \eqref{22} can be evaluated as
\begin{equation}
S_{n,l,m}^{(\al)}(x)=
\frac{(2\al+1)_{l+n} (2\al+1)_{m+n} (\al+\half)_{l+m}}
{2^{2n} (\al+\half)_l (\al+\half)_m (2\al+1)_{l+m} (\al+1)_n^2}\,(x^2-1)^n\,
R_{l-n}^{(\al+n,\al+n)}(x)\,R_{m-n}^{(\al+n,\al+n)}(x).
\label{45}
\end{equation}
\end{theorem}
\Proof
By \eqref{22}, \eqref{21} and \eqref{23} we obtain the recurrence
\begin{equation*}
S_{n,l,m}^{(\al)}
=\frac{-l-m-\al+\half}{\al+1}\,(1-x^2)\,S_{n-1,l-1,m-1}^{(\al+1)}.
\end{equation*}
Iteration gives
\begin{equation*}
S_{n,l,m}^{(\al)}=\frac{(-l-m-\al+\half)_n}{(\al+1)_n}\,(1-x^2)^n\,
S_{0,l-n,m-n}^{(\al+n)}.
\end{equation*}
Now use \eqref{22}, \eqref{18} and \eqref{29}.\qed
\bPP
As an immediate corollary, by the orthogonality relation 
\eqref{30} for Racah polynomials and by substitution of
\eqref{31} and \eqref{29}, we obtain:
\begin{theorem}[Dual addition formula]
For $l\ge m$ and for $j\in\{0,\ldots,m\}$ there is the expansion
\begin{multline}
R_{l+m-2j}^{(\al,\al)}(x)
=\sum_{n=0}^m\frac{\al+n}{\al+\thalf n}\,
\frac{(-l)_n (-m)_n (2\al+1)_n}{2^{2n}(\al+1)_n^2\,n!}\,
(x^2-1)^n\,
R_{l-n}^{(\al+n,\al+n)}(x)\,R_{m-n}^{(\al+n,\al+n)}(x)\\
\times R_n\big(j(j-l-m-\al-\thalf);\,
\al-\thalf,\al-\thalf,-m-1,-l-\al-\thalf\big).
\label{40}
\end{multline}
\end{theorem}

In particular, for $j=0$,
\begin{equation}
R_{l+m}^{(\al,\al)}(x)
=\sum_{n=0}^m
\frac{\al+n}{\al+\thalf n}\,
\frac{(-l)_n (-m)_n (2\al+1)_n}{2^{2n}(\al+1)_n^2\,n!}\,
(x^2-1)^n\,
R_{l-n}^{(\al+n,\al+n)}(x)\,R_{m-n}^{(\al+n,\al+n)}(x),
\end{equation}
and for $j=m$ we obtain by \eqref{25} that
\begin{equation}
R_{l-m}^{(\al,\al)}(x)
=\sum_{n=0}^m\binom mn
\frac{\al+n}{\al+\thalf n}\,
\frac{(l+2\al+1)_n (2\al+1)_n}{2^{2n}(\al+1)_n^2}\,
(1-x^2)^n\,
R_{l-n}^{(\al+n,\al+n)}(x)\,R_{m-n}^{(\al+n,\al+n)}(x),
\end{equation}
which is dual to \eqref{44} and
which has a further specialization to
\begin{equation}
1=\sum_{n=0}^m\binom mn
\frac{\al+n}{\al+\thalf n}\,
\frac{(m+2\al+1)_n (2\al+1)_n}{2^{2n}(\al+1)_n^2}\,
(1-x^2)^n\,
\big(R_{m-n}^{(\al+n,\al+n)}(x)\big)^2.
\label{43}
\end{equation}
Note that \eqref{43} coincides with formula \eqref{49}, which is
a specialisation of the addition formula \eqref{42} for ultraspherical
polynomials.
\begin{remark}
It follows from \eqref{22} and \eqref{57} that
\begin{multline}
\int_{-1}^1 S_{n,l,m}^{(\al)}(x)\,R_{l+m-2j}^{(\al,\al)}(x)\,
(1-x^2)^\al\,dx
=w_{\al-\half,\al-\half,-m-1,-l-\al-\half}(j)\,
h_{l+m-2j}^{(\al,\al)}\\
\times R_n\big(j(j-l-m-\al-\thalf);\,
\al-\thalf,\al-\thalf,-m-1,-l-\al-\thalf\big).
\label{58}
\end{multline}
By \eqref{45} and \eqref{59} we can rewrite this as
\begin{align*}
&\int_{-1}^1 R_{l-n}^{(\al+n,\al+n)}(x)\,R_{m-n}^{(\al+n,\al+n)}(x)\,
R_{l+m-2j}^{(\al,\al)}(x)\,(1-x^2)^{\al+n}\,dx\\
&\qquad\qquad\qquad\qquad=\const
\hyp43{-n,n+2\al,-j,j-l-m-\al-\thalf}{\al+\thalf,-l,-m}1\\
&\qquad\qquad\qquad\qquad=\const
\hyp43{-m+n,-m-n-2\al,j-m,l-j+\al+\thalf}
{-m,-m-\al+\thalf,l-m+1}1,
\end{align*}
where the second equality follows by twofold application of
Whipple's identity \cite[Theorem 3.3.3]{12} and
where the constants can be given as explicit, elementary, but somewhat
tedious expressions.
It turns out that the second ${}_4F_3$ evalutaion of the integral above
precisely matches the formula given by Koelink et al.\ \cite[(2.6]{17}.
Just put there (without loss of generality) $k=0$ and replace
$\al,\be,n,m,t$ by
$\al+n-\thalf,\al+\thalf,l-n,m-n,j-n$, respectively.
So, in a sense, the dual addition formula for ultraspherical polynomials
was already derived there in disguised form.
\end{remark}
\section{A limit to Hermite polynomials}
\label{62}
We will do a rescaling in the dual addition formula \eqref{40} such
that we can take the limit for $\al\to\iy$. For this purpose observe that
the Racah polynomial \eqref{59}
(where $l\ge m\ge \max(j,n)$) has limits
\begin{equation}
\begin{split}
\lim_{\al\to\iy}\al^{-j} R_n\big(j(j-l-m-\al-\thalf);\,
\al-\thalf,\al-\thalf,-m-1,-l-\al-\thalf\big)
&=\frac{2^j (-n)_j}{(-l)_j(-m)_j}\,,\\
\lim_{\al\to\iy}\al^{-n} R_n\big(j(j-l-m-\al-\thalf);\,
\al-\thalf,\al-\thalf,-m-1,-l-\al-\thalf\big)
&=\frac{2^n(-j)_n}{(-l)_n(-m)_n}\,.
\end{split}
\label{54}
\end{equation}
Otherwise said,
$R_n\big(j(j-l-m-\al-\thalf);\,
\al-\thalf,\al-\thalf,-m-1,-l-\al-\thalf\big)=O(\al^{\min(n,j)})$
as $\al\to\iy$ with the order constant given in \eqref{54}.

Now, in \eqref{40}, replace $x$ by $\al^{-\half}x$, multiply both sides
by $\al^{\half(l+m-2j)}$ and let $\al\to\iy$.
By \eqref{52} and \eqref{54} we obtain
\begin{equation}
2^j (-l)_j(-m)_j H_{l+m-2j}(x)=\sum_{n=j}^m
\frac{(-n)_j}{n!}\,(-2)^n (-l)_n(-m)_n\,
H_{l-n}(x) H_{m-n}(x)\quad(l\ge m),
\label{46}
\end{equation}
which may be called the dual addition formula for Hermite polynomials.
Formula \eqref{46} for arbitrary $j$ is equivalent to its case $j=0$,
\begin{equation}
H_{l+m}(x)=\sum_{n=0}^m
\frac{(-2)^n (-l)_n(-m)_n}{n!}\,
H_{l-n}(x) H_{m-n}(x)\qquad(l\ge m),
\end{equation}
and this is precisely \cite[10.13(36)]{16}.

Next we want to consider the limit as $\al\to\iy$ of \eqref{45} with
$S_{n,l,m}^{(\al)}$ given by \eqref{22}. Recall that \eqref{45} together with
\eqref{22} is
the dual of \eqref{40} in the sense of Fourier--Racah inversion.
Observe from \eqref{26}--\eqref{29} that
\begin{align}
\lim_{\al\to\iy} \al^j w_{\al-\half,\al-\half,-m-1,-l-\al-\half}(j)
&=\frac{(-l)_j (-m)_j}{2^j j!},\label{55}\\
\lim_{\al\to\iy} \al^{-n} h_{n;\al-\half,\al-\half,-m-1,-l-\al-\half}
&=\frac{2^n n!}{(-l)_n(-m)_n}\,.\label{56}
\end{align}
In \eqref{45}, replace $x$ by $\al^{-\half}x$, multiply both sides
by $\al^{\half(l+m-2n)}$ and let $\al\to\iy$.
By \eqref{52}, \eqref{54} and \eqref{55}
we obtain
\begin{equation}
\sum_{j=n}^m
\frac{(-j)_n}{j!}\,2^j (-l)_j(-m)_j \,
H_{l+m-2j}(x)=(-2)^n(-l)_n(-m)_n\,H_{l-n}(x) H_{m-n}(x)\quad(l\ge m).
\label{47}
\end{equation}
Again, as with \eqref{46}, formula \eqref{47} for arbitrary $n$
is equivalent to its case $n=0$,
\begin{equation}
\sum_{j=0}^m\frac{2^j (-l)_j(-m)_j}{j!}\,H_{l+m-2j}(x)=H_l(x) H_m(x)
\qquad(l\ge m),
\end{equation}
and this is precisely the linearization formula \cite[p.42]{10} for
Hermite polynomials.

Just as with \eqref{40} and \eqref{45}, the identities  \eqref{46} and
\eqref{47} can be obtained from each other by a Fourier type inversion.
This no longer involves an orthogonal system as the Racah polynomials but
a biorthogonal system implied by
the \emph{ biorthogonality relation}
(see Riordan \cite[Section 2.1]{15} or
Krattenthaler \cite[(1.1)]{24} with $a_j=1$, $b_j=0$)
\begin{equation}
\sum_{j=0}^\iy \frac{(-n)_j}{j!}\, \frac{(-j)_k}{k!}=\de_{n,k}\,.
\label{48}
\end{equation}
Note that the above sum in fact runs form $j=k$ to $n$.
For $k<n$ formula \eqref{48} is equivalent to
${}_1F_0(-n+k;-;1)=\sum_{j=0}^{n-k}\binom{n-k}j (-1)^j=0$.
\\[\smallskipamount]\indent
The biorthogonality \eqref{48} is also a limit case of
the Racah orthogonality relation \eqref{30}. Indeed, replace
$\al,\be,\ga,\de$ by $\al-\thalf,\al-\thalf,-m-1,-l-\al-\thalf$,
multiply both sides of \eqref{30} by $\al^{-n}$, let $\al\to\iy$,
and use \eqref{54}, \eqref{55} and \eqref{56}.
It is quite remarkable that a biorthogonal (and essentially non-orthogonal)
system can be obtained as a limit case of an orthogonal system.
Of course, before the limit it taken, the orthogonal system already
has to be prepared as a biorthogonal system by rescaling.
\section{The dual addition formula for Gegenbauer functions}
The dual product formula for the functions \eqref{1} is given in
\cite[(4.17)]{1} as
\begin{equation}
F(g;r,2p) F(g;r,2q)=\frac1{8\pi}\int_0^\iy F(g;r,2k)
\frac{\prod_{\de_1,\de_2,\de_3=+,-}
\Ga(\thalf(g+i\de_1 p+i\de_2 q+i\de_3 k))}
{\Ga(g)^2 \prod_{\de=+,-}\Ga(i\de k)\Ga(g+i\de k)}\,dk,
\label{6}
\end{equation}
where $g\in(0,\iy)$ and $r,p,q\in\RR$.
The formula is obtained there as a limit case of a similar formula
for a $q$-analogue (or relativistic analogue) of the Gegenbauer
function.
By \eqref{4} we can rewrite \eqref{6} as
\begin{multline}
\frac{\Ga(\al+\thalf)^2\,|\Ga(\al+\thalf+2i\la)|^2\,
|\Ga(\al+\thalf+2i\mu)|^2}{\Ga(2\al+1)}\,
\phi_{2\la}^{(\al,-\half)}(t)\,\phi_{2\mu}^{(\al,-\half)}(t)\\
=\frac1{4\pi}\int_{-\iy}^\iy \phi_{2\nu}^{(\al,-\half)}(t)\,
\Bigg|\frac{\Ga\big(i\nu\pm i\la\pm i\mu+\thalf\al+\tfrac14\big)}
{\Ga(2i\nu)}\Bigg|^2 d\nu,
\label{7}
\end{multline}
where $\al>-\thalf$ and $t,\la,\mu\in\RR$.
Note that the integral in \eqref{7} converges absolutely by
\eqref{32} and by estimates for the Gamma quotient using
\dlmf{5.5.5}, \dlmf{5.11.12} and, from \dlmf{5.5.3},
\[
|\Gamma(\thalf+i\nu)|^2=\frac\pi{\cosh(\pi\nu)}\,.
\]
The cases $\al=0$ and $\half$ of \eqref{7} were earlier given
by Mizony \cite{5}.

We recognize the weight function in the integrand of \eqref{7}
as the weight function in the orthogonality relation \eqref{8}
for Wilson polynomials with parameters 
$\pm i\la\pm i\mu+\thalf\al+\tfrac14$.
The case $t=0$ of \eqref{7} coincides with the case $m=n=0$ of \eqref{8}.

Similarly as the sum \eqref{22} is suggested by formula
\eqref{17}, we are led by formula \eqref{7} to try to evaluate
the integral
\begin{equation}
I_n^\al(\la,\mu):=
\frac1{4\pi}\int_{-\iy}^\iy \phi_{2\nu}^{(\al,-\half)}(t)\,
W_n(\nu^2;\pm i\la\pm i\mu+\thalf\al+\tfrac14)\,
\Bigg|\frac{\Ga\big(i\nu\pm i\la\pm i\mu+\thalf\al+\tfrac14\big)}
{\Ga(2i\nu)}\Bigg|^2 d\nu.
\label{9}
\end{equation}
Here and further in this section we will only work formally. We
will not bother about convergence, moving of integration contours
and justification of Fourier--Wilson inversion. Certainly this should
be repaired later.

By \eqref{33} and by shifting integration contours we get
\begin{align*}
I_n^\al(\la,\mu)&=
\frac1{4\pi}\int_{-\iy}^\iy 
\frac{\phi_{2\nu-i}^{(\al,-\half)}(t)-\phi_{2\nu+i}^{(\al,-\half)}(t)}
{2i\nu}\,
W_{n-1}(\nu^2;\pm i\la\pm i\mu+\thalf(\al+1)+\tfrac14)
\\
&\qquad\qquad\qquad\qquad
\times\Bigg|\frac{\Ga\big(i\nu\pm i\la\pm i\mu+\thalf(\al+1)+\tfrac14\big)}
{\Ga(2i\nu)}\Bigg|^2 d\nu
=\frac{\sinh^2 t}{\al+1}\,I_{n-1}^{\al+1}(\la,\mu),
\end{align*}
where the last equality follows by \eqref{34}.
Iteration gives
\begin{equation*}
I_n^\al(\la,\mu)=\frac{(\sinh t)^{2n}}{(\al+1)_n}\,I_0^{\al+n}(\la,\mu).
\end{equation*}
Hence, by \eqref{9} and \eqref{7},
\begin{multline}
\frac1{4\pi}\int_{-\iy}^\iy \phi_{2\nu}^{(\al,-\half)}(t)\,
W_n(\nu^2;\pm i\la\pm i\mu+\thalf\al+\tfrac14)\,
\Bigg|\frac{\Ga\big(i\nu\pm i\la\pm i\mu+\thalf\al+\tfrac14\big)}
{\Ga(2i\nu)}\Bigg|^2 d\nu\\
=\frac{\Ga(\al+\thalf)^2\,|\Ga(n+\al+\thalf+2i\la)|^2\,
|\Ga(n+\al+\thalf+2i\mu)|^2}{\Ga(2n+2\al+1)}\,
\frac{(\sinh t)^{2n}}{(\al+1)_n}\,
\phi_{2\la}^{(\al+n,-\half)}(t)\,\phi_{2\mu}^{(\al+n,-\half)}(t).
\label{13}
\end{multline}
By \eqref{8} there corresponds, at least formally, to \eqref{13}
the orthogonal expansion
\begin{equation*}
\phi_{2\nu}^{(\al,-\half)}(t)=\sum_{k=0}^\iy
\frac{(\sinh t)^{2k}}{(\al+1)_k (k+2\al)_k\,k!}
\phi_{2\la}^{(\al+k,-\half)}(t)\,\phi_{2\mu}^{(\al+k,-\half)}(t)\,
W_k(\nu^2;\pm i\la\pm i\mu+\thalf\al+\tfrac14).
\end{equation*}
Equivalently, by \eqref{16}, we can write what we call the
\emph{ dual addition formula for Gegenbauer functions}:
\begin{equation}
\phi_{4\nu}^{(\al,\al)}(t)=\sum_{k=0}^\iy
\frac{(\sinh 2t)^{2k}}{(\al+1)_k (k+2\al)_k\,k!}\,
\phi_{4\la}^{(\al+k,\al_+k)}(t)\,\phi_{4\mu}^{(\al+k,\al+k)}(t)
\,W_k(\nu^2;\pm i\la\pm i\mu+\thalf\al+\tfrac14).
\label{15}
\end{equation}
\section{Further perspective}
\label{64}
The results of this paper suggest much further work.
I will only discuss here the polynomial case.
A very obvious thing to do is to imitate the approach of Section 4
for $q$-ultraspherical polynomials, starting with their
linearization formula \cite[(10.11.10)]{12}, which goes back to
Rogers (1895). This was settled by the author in \cite{25}.
The $q$-Racah polynomials pop up there. Furthermore, this dual addition formula
can also be proved from the Rahman--Verma
addition formula \cite{13} for $q$-ultraspherical
polynomials.

It would be very interesting to find a dual addition formula
for the addition formula for Jacobi polynomials \cite{22},
starting with the linearization formula \cite{23}.

It would be quite challenging to search for an addition formula on a
higher
level (for \mbox{($q$-)}Racah polynomials) which has both the addition formula
and the dual addition formula for ultraspherical polynomials as a limit case.
The recent formula for the linearization coefficients of
Askey--Wilson polynomials by Foupouagnigni et al.\
\cite[Theorem 21]{14}, which involves four summations, does not give much hope
for a quick answer to this problem.

Very important will also be to give a group theoretic interpretation
for the dual addition formula for ultraspherical polynomials, for instance
when $\al=0$. Possibly this can be done in the context of
tensor algebras associated with the group SU(2).
\subsection*{Acknowledgements}
The author wants to thank Erik Koelink (Nijmegen)
for pointing out that
\eqref{45} would possibly be equivalent to \cite[(2.6]{17}
and for suggesting to consider the Hermite limit of the dual
addition formula. Josef Hofbauer (Vienna) was very helpful
in providing material about the history of the addition formula
for ultraspherical polynomials.

\quad\\
\begin{footnotesize}
\begin{quote}
{ T. H. Koornwinder, Korteweg-de Vries Institute, University of
 Amsterdam,\\
 P.O.\ Box 94248, 1090 GE Amsterdam, The Netherlands;

\vspace{\smallskipamount}
email: }\texttt{thkmath@xs4all.nl}
\end{quote}
\end{footnotesize}

\end{document}